\documentclass[12pt]{article}

\usepackage{url}
\usepackage[normalem]{ulem}
\usepackage{amsmath}
\newcommand{\stkout}[1]{\ifmmode\text{\sout{\ensuremath{#1}}}\else\sout{#1}\fi}\usepackage{amssymb}
\usepackage{amscd}
\usepackage{tikz}
\usetikzlibrary{arrows,automata,positioning}

\usepackage{amsmath}
\usepackage{amsthm}
\usepackage{amsfonts}
\usepackage{bm}

\newtheorem{Theorem}{Theorem}
\newtheorem{lemma}[Theorem]{Lemma}
\newtheorem{corollary}[Theorem]{Corollary}

\theoremstyle{remark}
\newtheorem{remark}[Theorem]{Remark}

\numberwithin{Theorem}{section}

\begin{document}
\title{The analogue of overlap-freeness for the period-doubling sequence}
\author{James D. Currie\footnote{Supported by the Natural Sciences and Engineering Research Council of Canada (NSERC), [funding
reference number 2017-03901].}\\
Department of Mathematics and Statistics\\
University of Winnipeg\\
Winnipeg, Manitoba R3B 2E9, Canada}
\maketitle
\begin{abstract}
Good words are binary words avoiding factors 11 and 1001, and patterns 0000 and 00010100. We show that good words bear the same relationship to the period-doubling sequence that overlap-free words bear to the Thue-Morse sequence. We prove an analogue of Fife's Theorem for good words, exhibit the lexicographically least and greatest infinite good words, and determine the patterns avoided by the period doubling word.
\vspace{.1in}

\end{abstract}

\section{Introduction}
The famous Thue-Morse sequence ${\bf t}$ is a fixed point of the binary morphism $\mu$ given by $\mu(0)=01$, $\mu(1)=10$, namely 
\[{\bf t}=\lim_{n\rightarrow\infty}\mu^n(0)\]

Thue \cite{thue06} proved that ${\bf t}$ is overlap-free. 
\begin{Theorem}\label{thue freeness} Let $w$ be an overlap-free binary word. Then $\mu(w)$ is overlap-free.
\end{Theorem}
He also showed that, in the case of two-sided infinite words and circular words, every overlap-free binary word is the image under $\mu$ of an overlap-free word \cite{thue12}. The analysis of words with `ends' is more complicated, but it has long been known that finite overlap-free binary words also arise via iterating $\mu$. (See Restivo and Salemi \cite{restivo84} for example.)
\begin{Theorem}\label{thue factorization} Let $w\in\{0,1\}^*$ be a finite overlap-free word. Then we can write $w=a\mu(u)b$, where $a, b\in\{\epsilon,0,00,1,11\}$, and $u$ is overlap-free. If $|w|\ge 7$ this factorization is unique. If ${\bm w}$ is a one-sided infinite overlap-free word, then we can write ${\bm w}=a\mu({\bm u})$, for some one-sided infinite overlap-free word ${\bm u}$ where $a\in\{\epsilon,0,00,1,11\}$.
\end{Theorem}

Restivo and Salemi used their version of Theorem~\ref{thue factorization} to give a rough  enumeration of binary overlap-free words. A better enumeration was given by Kobayashi \cite{kobayashi88}, who got a good lower bound by counting finite words which extend to infinite overlap-free words. For this, Kobayashi used the deep theorem of Fife  \cite{fife80} characterizing the infinite overlap-free words.  Carpi \cite{carpi93} and Cassaigne \cite{cassaigne93} extended Fife's theorem by characterizing
the set of finite overlap-free binary words via a regular language. 

Because of Theorem~\ref{thue factorization}, ${\bf t}$ turns up frequently in the study of overlap-free binary words. An example is the following result of Berstel \cite{berstel94} (later greatly generalized by Allouche {\em et al.} \cite{allouche98}):

\begin{Theorem}\label{lex}
The lexicographically greatest one-sided infinite overlap-free binary word starting with 0 is 
${\bf t}$.
\end{Theorem}

The relationship between ${\bf t}$ and $\mu$ was also key to establishing the following theorem of Shur \cite{shur96}:

\begin{Theorem} Suppose ${\bf t}$ encounters a pattern $p\in\{0,1\}^*$. Then either $p$ is a factor of ${\bf t}$, or $p$ is one of 00100 and 11011. 
\end{Theorem}

Another famous binary sequence is the period-doubling sequence, which is the fixed point
\[{\bm d}=01000101010001000100010101000101\cdots\]
of the binary morphism $\delta$ where $\delta(0)=01$, $\delta(1)=00$.  This sequence has been much studied in the context of quasicrystal spectral theory. (See Damanik \cite{damanik00}, for example.)

Call a binary word $w$ {\bf good} if it does not contain the factor 11 or 1001, and does not encounter either of the patterns 0000 or 00010100. We will show that ${\bm d}$ is good. In fact, we show that the period doubling morphism $\delta$ has the same relationship to good words that $\mu$ has to overlap-free words, namely:

\begin{Theorem}\label{f good}Suppose $w$ is good. Then $\delta(w)$ is good.
\end{Theorem}
\begin{Theorem}\label{main theorem} Let $w$ be a finite good word. Then we can write $w=a\delta(u)b$, where $a\in\{\epsilon,0,1\}$, $b\in\{\epsilon,0\}$, and $u$ is good. If $|w|\ge 4$ this factorization is unique. If ${\bm w}$ is a one-sided infinite good word, then we can write ${\bm w}=a\delta({\bm u})$, some one-sided infinite good word ${\bm u}$ where $a\in\{\epsilon,0,1\}$.
\end{Theorem}

We build on these theorems to:
\begin{itemize}
\item Give a version of Fife's theorem for good words, characterizing infinite good words;
\item Exhibit lexicographically extremal one-sided infinite good words; 
\item Characterize the binary patterns avoided by ${\bm d}$.
\end{itemize}
\section{Good words}

Unless otherwise specified, our words and morphisms are over the binary alphabet $\{0,1\}$. It is convenient to record morphisms inline, i.e., $g=[g(0),g(1)]$.

\begin{lemma}\label{delta^-1 good}
Let $u$ be a finite binary word. Suppose $\delta(u)$ is good. Then $u$ is good. \end{lemma}
\begin{proof} If $u$ encounters pattern 0000 or 00010100, so does $\delta(u)$. If $u$  has factor 11 or 1001, then $\delta(u)$ has factor $\delta(11)=0000$ or $\delta(1001)=00010100$, and thus encounters pattern 0000 or 00010100.\end{proof}

\begin{remark} \label{parity} Let $w=w_1w_2w_3\cdots w_{2n}$ with $w_i\in\{0,1\}$. Word $w$ can be written as $w=\delta(u)$ for some $u$ if and only if $w_i=0$ for each odd index $i$.
\end{remark}
\begin{lemma}\label{factorization}
Let $w$ be a finite binary word with no factor 11, 1001, or 0000. Then we can write $w=a\delta(u)b$ where $a\in\{\epsilon,0,1\}$ and $b\in\{\epsilon,0\}.$ If $|w|\ge 4$ this factorization is unique. If ${\bm w}$ is a one-sided infinite word with no factor 11, 1001, or 0000, then we can write ${\bm w}=a\delta({\bm u})$, some one-sided infinite word ${\bm u}$ where $a\in\{\epsilon,0,1\}$.
\end{lemma}
\begin{proof} First we demonstate the existence of the factorization for finite words:  If $|w|_1\le 1$ then $w$ is a factor of 0001000, and the result is established by a finite check. Suppose $|w|_1\ge 2$. If $10^k1$ is a factor of $w$, the conditions on $w$ force $k=1$ or $k=3$. By induction, if $1u1$ is a factor of $w$ then $|1u|$ is even; therefore, all the 1's in $w$ have index in $w$ of the same parity. If the parity of the indices of 1's is even, let $|a|=0$; let $|b|$ be 0 (resp., 1) if $|w|$ is even (resp., odd).  If the parity of the indices of 1's is odd, let $|a|=1$; let $|b|$ be 1 (resp., 0) if $|w|$ is even (resp., odd). By Remark~\ref{parity} we can write $a^{-1}wb^{-1}=\delta(u)$ for some $u$. 

We have shown that we can write $w=a\delta(u)b$ and $a,b\in\{\epsilon, 0,1\}$. Suppose that $b=1$. Recall that $|w|_1\ge 2$. The parity of the indexes of $a$ and $b$ is different, so we cannot have $a=1$. It follows that $|\delta(u)|_1\ge 1$. Then $\delta(u)1$ must have suffix 011, 01001, or 00001, none of which are good. This is a contradiction so in fact $b\in \{\epsilon, 1\}$.
 
Now we show uniqueness of the factorization for finite words: Suppose $|w|\ge 4$ and $w$ has two factorizations $w=a_1\delta(u_1)b_1=a_2\delta(u_2)b_2$. If $|a_1|=|a_2|$ we are forced to choose $a_1=a_2$, $u_1=u_2$, and $b_1=b_2$, so the factorizations are identical. Suppose without loss of generality then that $|a_1|=0$, $|a_2|=1$. Then by Lemma~\ref{parity}, any 1's in $a_1\delta(u_1)b_1$ have even index, but 
any 1's in $a_2\delta(u_2)b_2$ have even index. Since $a_1\delta(u_1)b_1=w=a_2\delta(u_2)b_2$, we conclude that $|w|_1=0$, so that $w$ has 0000 as a prefix, which is impossible. 

Now suppose that ${\bm w}$ is a one-sided infinite good word. For each non-negative $n$, let the length $n$ prefix of ${\bm w}$ be $w_n$. 
We have proved that finite words can be factored, so write each $w_n=a_n\delta(u_n)b_n$ where $a_n,b_n\in\{\epsilon,0,1\}.$ Word $w_4$ cannot be 0000, so that $|w_4|_1>0$. By Remark~\ref{parity}, the index of the first 1 in $w_4$ determines $a_4$ and all subsequent $a_n$, so that $a_n=a_4$ for $n\ge 4$. This implies that $\delta(u_n)$ is a prefix of $\delta(u_{n+1})$ for $n\ge 4$, so that $u_{n+4}$ is a prefix of $u_{n+5}$ for all $n$. 
Let ${\bm u}=\lim_{n\rightarrow\infty} u_{n+4}$.
Then
\begin{eqnarray*}{\bm w}&=&\lim_{n\rightarrow\infty} w_n\\
&=&\lim_{n\rightarrow\infty} w_{n+4}\\
&=&\lim_{n\rightarrow\infty} w_{n+4}b_{n+4}^{-1}\\
&=&\lim_{n\rightarrow\infty} a_4\delta(u_{n+4})\\
&=&a_4\lim_{n\rightarrow\infty} \delta(u_{n+4})\\
&=&a_4\delta(\lim_{n\rightarrow\infty} u_{n+4})\\
&=&a_4\delta({\bm u}).
\end{eqnarray*}
\end{proof}

\noindent{\bf Proof of Theorem~\ref{main theorem}:} This is immediate from Lemma~\ref{delta^-1 good} and Lemma~\ref{factorization}.$\qed$
\vspace{.05in} 

Call a non-erasing morphism $g$ {\bf even} if, for every letter $u$, $|g(u)|$ is even.
\begin{lemma}\label{even lengths}
Let $p$ be a  pattern and let $w$ be a binary word. Suppose that $g(p)$ is a factor of $\delta(w)$ where $g$ is an even morphism. Then $w$ encounters pattern $p$.
\end{lemma}
\begin{proof}Write $p=u_1u_2\cdots u_n$ where the $u_i$ are letters. Write $\delta(w)=aU_1U_2\cdots U_nb$, where $U_i=g(u_i)$ for each $i$. If $|a|$ is even, then $w=\delta^{-1}(a)\delta^{-1}(U_1U_2\cdots U_n)\delta^{-1}(b)$, and $w$ contains the instance $\delta^{-1}(g(p))$ of $p$.

If $|a|$ is odd, then 0 is the last letter of $a$ and of each $U_i$. Thus
$\delta(w)=a0^{-1}0U_10^{-1}0U_20^{-1}0\cdots U_n0^{-1}0b$, 
and $w$ contains the instance $\delta^{-1}(h(p))$ of $p$, where $h$ is the morphism defined on the letters of $p$ by $h(u_i)=0g(u_i)0^{-1}$.\end{proof}
\begin{lemma}\label{uvu} Let $w,u,v$ be binary words and suppose that $uvu$ is a factor of $\delta(w)$. If $|u|_1>0$ then $|uv|$ is even.
\end{lemma}

\begin{proof} Since $uvu$ has period $|uv|$, $\delta(w)$ has a factor $1z1$ where $|1z|=|uv|$. The Lemma follows by Remark~\ref{parity}.
\end{proof}

\noindent {\bf Proof of Theorem~\ref{f good}:} To begin with, we show that $\delta(w)$ does not contain any of 11, 0000, 1001, or 00010100 as a factor. By Remark~\ref{parity}, $\delta(w)$ does not have a factor 11 or 1001. If $\delta(w)$ has factor 0000, write $\delta(w)=a0000b$ for words $a$ and $b$. If $|a|$ is even, then $w$ has prefix $\delta^{-1}(a0000)$, which ends in 11. This is impossible; if $|a|$ is odd, then the last letter of $a$ is 0, so that $w$ has prefix $\delta^{-1}((a0^{-1})0000)$, which again ends in 11. Finally, if 00010100 is a factor of $\delta(w)$, then $w$ contains factor $\delta^{-1}(00010100)=1001$, which is impossible. (The index of 00010100 in $\delta(w)$ must be odd by Remark~\ref{parity})

Suppose now that $\delta(w)$ encounters pattern $p=0000$, so that $XXXX$ is a factor of $\delta(w)$ for some non-empty $X$. Since $\delta(w)$ doesn't have 0000 as a factor, we must have $|X|_1>0$. Using $u=X$ and $v=\epsilon$ in Lemma~ \ref{uvu}, we conclude that $|X|$ is even. Then Lemma~ \ref{even lengths}  implies that $w$ encounters $0000$, which is a contradiction.

Suppose that $\delta(w)$ encounters pattern $p=00010100$, so that $XXXYXYXX$ is a factor of $\delta(w)$ for some non-empty $X$ and $Y$. Suppose that $|X|_1>0$. Since $XX$ is a factor of $\delta(w)$,  letting $u=X$ and $v=\epsilon$ in Lemma~ \ref{uvu} implies that $|X|$ is even. Again, since $XYX$ is a factor of $\delta(w)$,  letting $u=X$ and $v=Y$ in Lemma~ \ref{uvu} implies that $|XY|$ is even. It follows that $|Y|$ is even. Then Lemma~ \ref{even lengths}  implies that $w$ encounters $00010100$, which is a contradiction. We therefore conclude that $|X|_1=0$, so that $X=0^n$ for some $n\ge 1$. Since $XX=0^{2n}$ is a factor of $\delta(w)$, but 0000 is not, we conclude that $n=1$ and $X=0$. Thus $\delta(w)$ contains the factor $XXXYXYXX=000Y0Y00$.

Since 0000 is not a factor of $\delta(w)$, the first letter of $Y$ is 1. If $Y=1$, then 
$\delta(w)$ has factor 00010100, which is impossible. Therefore $|Y|>1.$

Suppose $Y$ ends in 0. If $Y$ ends in 10, then $\delta(w)$ contains factor 1001 (inside $Y0Y$), which is impossible. If $|Y|$ ends in 00, then $000Y0Y00$ ends in $Y00$, hence 0000, which is again impossible. Thus $Y$ ends in 1, hence in 01. Write $Y=1Z1$ for some non-empty word $Z$. Word $\delta(w)$ has the factor $0001Z101Z100$.

If $|1Z1|=3$, then $\delta(w)$ contains $01Z101Z1=01010101$, an instance of $0000$, already proved impossible. It follows that $|1Z1|\ge 4$. If $Z$ begins 01, then $Z10$ begins either 0101 or 0100. However, if $Z10$ begins 0101 then, since $Z$ ends in 0, $Z101Z10$ again contains 01010101; if $Z10$ begins 0100 then $0001Z10$ begins 00010100, which is not a factor of $\delta(w)$.
We conclude that $Z$ does not begin 01 and therefore begins 00. 

If $Z10$ begins 001 then 1Z10 has the impossible factor 1001. Thus  
$Z10$ begins 0001. We now consider suffixes of $Z.$ 
If $Z$ ends 10, then $01Z$ ends either 1010 or 0010. However, if $01Z$ ends 1010 then, since $Z$ begins with 0, $01Z101Z$ contains 10101010, an instance of $0000$, which is impossible; if $01Z$ ends in 0010 then $101Z100$ ends 00010100, which is not a factor of $\delta(w)$.
We conclude that $Z$ doesn't end 10 and therefore ends 00.
If $01Z$ ends 0100 then $\delta(w)$ contains $01Z1$, hence the impossible factor 1001. Thus  $01Z$ ends 1000. Since $Z$ begins 00, this forces $01Z101Z$ to contain 00010100, which is impossible.\hfill$\qed$

\begin{corollary} The period-doubling word ${\bm d}$ is good.
\end{corollary}
\begin{Theorem}\label{o.s.i. good} Let ${\bf w}$ be a one-sided infinite binary word. Then ${\bm w}$ is good if and only if $\delta({\bm w})$ is good.
\end{Theorem}
\begin{proof}
 If $w_n$ is the length $n$ prefix of ${\bf w}$ then ${\bf w}=\lim_{n\rightarrow\infty}w_n$ and 
$
\delta({\bf w})=\lim_{n\rightarrow\infty}\delta(w_n).$ 
By Theorem~\ref{f good} and Lemma~\ref{delta^-1 good}, $w_n$ is good if and only if $\delta(w_n)$ is good.
\begin{eqnarray*}
{\bf w}\mbox{ is good }&\iff&\mbox{each }w_n\mbox{ is good }\\
&\iff&\mbox{each }\delta(w_n)\mbox{ is good }\\
&\iff&\lim_{n\rightarrow\infty}\delta(w_n)=\delta(\lim_{n\rightarrow\infty}w_n)=\delta({\bf w})\mbox{ is good.}\\
\end{eqnarray*}
\end{proof}

\begin{remark}\label{symmetry} While the set of binary overlap-free words is closed under complementation and reversal, the same is not true of good words. For example, 00101000 is good, but neither its complement nor its reversal is good.
\end{remark}

\section{An analogue of Fife's Theorem}

We characterize the one-sided infinite good words, developing a theory analogous to Fife's \cite {fife80} theory for  overlap-free binary words; however, we follow Ramersad's \cite{rampersad07} exposition of Fife, with appropriate modifications, rather than Fife's original paper. Let $G$ be the set of one-sided infinite good words.


\begin{lemma}\label{slide}
Suppose ${\bm u}$ is a one-sided infinite good word. Then $1{\bm u}$ is good if and only if $01{\bm u}$ is good.
\end{lemma}
\begin{proof}
Clearly if $01{\bm u}$ is good then $1{\bm u}$ is good. Suppose that $1{\bm u}$ is good but $01{\bm u}$ is  not. By Theorem~\ref{main theorem}, write $1{\bm u}=1\delta({\bm v})$ for some ${\bm v}$.
Since $1{\bm u}$ is good, $01{\bm u}$ must have a prefix which is either factor 11 or 1001, or a pattern instance $g(0000)$ or $g(00010100)$ where $g=[X,Y]$ is some non-erasing morphism. Clearly neither of 11 and 1001 can be a prefix of $01{\bm u}$, so $01{\bm u}$ has a prefix $XXXX$ or $XXXYXYXX$ where $X$ and $Y$ are non-empty. 

Suppose $01{\bm u}$ has prefix $XXXX$. Write $X=0X'$. Since $|XXXX|$ is even, word $01{\bm u}=\delta(0{\bm v})$ has $XXXX0=0X'0X'0X'0X'0$ as a prefix. But now the good word  
$1{\bm u}$ contains the fourth power $X'0X'0X'0X'0$, which is a contradiction.

Now  suppose that $01{\bm u}$ has prefix $XXXYXYXX$. Write $X=0X'$. If $X'=\epsilon$, then $X=0$. However then the length 2 prefix of  $01{\bm u}$ is $XX=00$, which is impossible. Thus $X'\ne\epsilon$, forcing 01 to be a prefix of $X$, so that $|X|_1>0$. Since $01{\bm u}=\delta(0{\bm v})$, considering the second $XX$ in $XXXYXYXX$, by Lemma~\ref{uvu} with $u=X$, $v=\epsilon$ we find that $|X|$ is even. Again,  $XYX$ is a factor of $\delta(0{\bm v})$, so applying Lemma~\ref{uvu} with $u=X$, $v=Y$ shows that $|XY|$ is also even. Since both $|X|$ and $|XY|$ are even, $|Y|$ is even. 
Since $|XXX|$ is even, $XXX0$ is a prefix of $01{\bm u}=\delta(0{\bm v})$, so that we can write 
$Y=0Y'$. Since $|XXXYXYXX|$ is even, word $01{\bm u}=\delta(0{\bm v})$ has $XXXYXYXX0=0X'0X'0X'0Y'0X'0Y'0X'0X'0$ as a prefix. But now the good word  
$1{\bm u}$ contains $X'0X'0X'0Y'0X'0Y'0X'0X'0=g'(00010100)$ where $g'=[X'0,Y'0]$. This is a contradiction.
\end{proof}
\begin{remark} This result uses the fact that  ${\bm u}$ is one-sided infinite. If $u=010101$, then $1u$ is good, but $01u=(01)^4$  is not.
\end{remark}

Let $G$ be the set of one-sided infinite good words. For $w\in\{0,1\}^*$, let $G_w=G\cap w\{0,1\}^\omega$.

\begin{lemma}\label{allouche}Let ${\bm w}$ be a one-sided infinite binary word.
\begin{enumerate}
\item[(a)] $\delta({\bm w})\in G\iff {\bm w}\in G$
\item[(b)] $1\delta({\bm w})\in G\iff 0{\bm w}\in G$
\item[(c)] $0\delta({\bm w})\in G\iff$ ($1{\bm w}\in G$) or (${\bm w}\in G_{001}$)
\end{enumerate}
\end{lemma}
\begin{remark} The cases in (c) are disjoint, since if 001 is a prefix of ${\bf w}$, then $1{\bf w}$ has prefix 1001 and is not good.
\end{remark}
\noindent {\bf Proof of (a):} This is just Theorem~\ref{o.s.i. good}. 

\noindent {\bf Proof of (b):} If $0{\bm w}\in G$, then by Theorem~\ref{o.s.i. good}, $\delta(0{\bm w})=01\delta({\bm w})\in G$, so in particular $1\delta({\bm w})\in G$. 

In the other direction, suppose $1\delta({\bm w})\in G$. Applying Lemma~\ref{slide} to prefixes, $01\delta({\bm w})=\delta(0{\bm w})\in G$. By Theorem~\ref{o.s.i. good}, $0{\bm w}\in G$.

\noindent {\bf Proof of (c):} If $1{\bm w}\in G$, then by Theorem~\ref{o.s.i. good} $\delta(1{\bm w})=00\delta({\bm w})\in G$, so in particular $0\delta({\bm w})\in G$. 

Suppose ${\bm w}\in G_{001}$. By Theorem~\ref{o.s.i. good}, $\delta({\bm w})\in G$. Suppose $0\delta({\bm w})\not\in G$. It therefore has a prefix $XXXX$ or $XXXYXYXX$ where $X$ is non-empty. Since 001 is a prefix of ${\bf w}$, $p=0010100$ is a prefix of $0\delta({\bf w})$. The prefix $XXX$ of $0\delta({\bf w})$ has period $|X|$. If $|X|<|p|$, then $|X|$ is a period of $p$, which has least period 5. We conclude that $|X|\ge 5$. This implies that, $|X|_1>0$. The second $XX$ in $XXX$ is a factor of $\delta({\bf w})$. By Lemma~\ref{uvu} with $u=X$, $v=\epsilon$, we conclude that $|X|$ is even. Therefore $|X|\ge 6$. Since $0^{-1}X$ is an odd-length prefix of $\delta({\bm w})$, 0 is a suffix of $X$. 

If $|X|=6$ then the prefix $0^{-1}XXX$ of $\delta({\bf w})$ ends in $0XX=0$ $001010$ $001010$, which has prefix 00010100. This is impossible, since $\delta({\bf w})\in G$. If $|X|>6$, then $p=0010100$ is a prefix of $X$, and 0 is a suffix of $X$, so factor $XX$ of $\delta({\bf w})$ contains $0$ $0010100$, again an impossibility.

In the other direction, suppose that $0\delta({\bf w})\in G$. Then ${\bm w}\in G$ by Theorem~\ref{o.s.i. good}, and we show that either $1{\bm w}\in G$, or  ${\bm w}\in 001(0,1)^*$. Suppose that $1{\bm w}\not\in G$ and 001 is not  a prefix of ${\bf w}$. It follows that 
a prefix of $1{\bm w}$ has the form 11, 1001, $XXXX$, or $XXXYXYXX$ with non-empty $X$ and/or $Y$. 

By Theorem~\ref{main theorem}, write ${\bm w}=a\delta({\bm u})$ where $a\in\{\epsilon,0,1\}$. 

If $11$ is a prefix of $1{\bm w}$, then $a=1$, and 10 is a prefix of  ${\bm w}$. However, then $0\delta({\bm w})$ has prefix $0\delta(10)=00001$, and is not good.

If $1001$ is a prefix of $1{\bm w}$, then  001 is a prefix of ${\bf w}$, which is a contradiction.

Suppose $XXXX$ is a non-empty prefix of $1{\bm w}$. Then the final $XX$ of $XXXX$ is a factor of $\delta({\bm u})$ and $|X|_1>0$, so by Lemma~\ref{uvu}, $|X|$ is even. Since $X$ is a prefix of $1{\bm w}$, the first letter of $X$ is 1. Since $XX$ is a factor of ${\bm w}$, and ${\bm w}$ has no factor 11, the last letter of $X$ is 0. Write $X=1X'0$.

Word ${\bm w}$ has prefix $X'01X'01X'01X'0c$ for some $c\in\{0,1\}$.
This implies that $0\delta({\bm w})$ has prefix
\begin{eqnarray*}
0\delta(X'01X'01X'01X'0c)&=&0\delta(X'0)\delta(1)\delta(X'0)\delta(1)\delta(X'0)\delta(1)\delta(X'0)\delta(c)\\
&=&0\delta(X'0)00\delta(X'0)00\delta(X'0)00\delta(X'0)0\overline{c}\\
&=&(0\delta(X'0)0)^4\overline{c}\\
\end{eqnarray*}
and $0\delta({\bm w})$ contains a fourth power. This is a contradiction, since $0\delta({\bm w})\in G$.

Suppose $XXXYXYXX$ is a prefix of $1{\bm w}$. This means that $XX$ and $XYX$ are factors of $\delta({\bm u})$, and $|X|_1>0$. We conclude by Lemma~\ref{uvu} that $|X|$ and $|Y|$ are even. The first letter of $X$ is 1. Since $XX$ and $YX$ are factors of ${\bm w}$, and ${\bm w}$ has no factor 11, the last letter of each of  $X$ and $Y$ is 0. Write $X=1X'0$ and $Y=cY'0$ where $c\in\{0,1\}$. Since $X'01$ is an even length prefix of $a\delta({\bm u})$, Remark~\ref{parity} forces $a=\epsilon$.

Then 
${\bm w}$ has prefix $X'01X'01X'01Y'01X'01y'01X'01X'0d$ for some $d\in\{0,1\}$.
Therefore $0\delta({\bm w})$ has prefix
\begin{eqnarray*}
&&0\delta(X'01X'01X'0cY'01X'0cY'01X'01X'0d)\\
&=&        0\delta(X')010
               0\delta(X')010
               0\delta(X')010
\overline{c}\delta(Y')010
               0\delta(X')010
\overline{c}\delta(Y')010
               0\delta(X')010
               0\delta(X')010\overline{d}\\
&=&h(00010100)\overline{d}
\end{eqnarray*}
where $h=[0\delta(X')010,\overline{c}\delta(Y')010]$. This is a contradiction, since $0\delta({\bm w})\in G$.$\qed$

Suppose $w\in\{0,1\}^*$ has a suffix $\delta^n(01)$, $n\ge 0$, and let $n$ be a as large as possible. Write $w=y\delta^{n}(01)$. Define mappings $\alpha$, $\beta$ and $\gamma$ on $w$ by
\begin{eqnarray*}
\alpha(w)&=&w\delta^n(00)=y\delta^{n+1}(01)\\
\beta(w)&=&w\delta^n(0100)=y\delta^{n+1}(001)\\
\gamma(w)&=&w\delta^n(010100)=y\delta^{n+1}(0001)\\
\end{eqnarray*}

For example, if $w=0001000101$, then $y=00$, $n=2$, $\delta^n(0)=0100$, $\delta^n(1)=0101$, so that 
\begin{eqnarray*}
\alpha(w)&=&00\mbox{ }0100\mbox{ }0101\mbox{ }0100\mbox{ }0100\\
\beta(w)&=&00\mbox{ }0100\mbox{ }0101\mbox{ }0100\mbox{ }0101\mbox{ } 0100\mbox{ }0100\\
\gamma(w)&=&00\mbox{ }0100\mbox{ }0101\mbox{ }0100\mbox{ }0101\mbox{ }0100\mbox{ }0101\mbox{ } 0100\mbox{ }0100\\
\end{eqnarray*}

Let ${\bm f}=f_1f_2f_3\cdots, \mbox{ where each } f_i\in B=\{\alpha,\beta,\gamma\}$. Suppose $w$ has some suffix $\delta^n(01)$. Then $w$ is a proper prefix of $f_1(w)$, which is a proper prefix of $f_2(f_1(w))$, which is a proper prefix of $f_3(f_2(f_1(w)))$, etc. We define the infinite composition
${\bm x} = (\cdots\circ f_3 \circ f_2 \circ f_1)(01)$ to be the one-sided infinite word
\[{\bm x}=\lim_{n\rightarrow\infty} f_n(f_{n-1}(\cdots f_2(f_1(w))\cdots )\] 
which has each $f_n(f_{n-1}(\cdots f_2(f_1(w))\cdots )$ as a prefix. Following Rampersad, we use the notation ${\bm x} = w\bullet{\bm f}$.

Let $w\in B^k$. Then 

\begin{equation}\label{2.3 01}
01\bullet w{\bm f} = (01\bullet w)\delta^k(01)^{-1}\delta^k({\bm x})\end{equation}

Define the sets $I$ and $F$ by
\[I =(\beta+\gamma)(\alpha\alpha)^*\alpha(\gamma+\beta\beta)\cup\gamma(\alpha\beta)^*\alpha\gamma \]
\[F = B^\omega-B^*IB^\omega\]
\begin{Theorem}\label{fife01}
Let ${\bm x}\in\{0, 1\}^\omega$.
If ${\bm x}$ begins with 01, then ${\bm x}$ is good if and only if ${\bm x}= 01\bullet {\bm f}$ for some
${\bm f}\in F$.
\end{Theorem}

We follow the notation of Berstel, also used by Rampersad. Here $I$ stands for `ideal', and $B^*IB^\omega$ is the ideal generated by $I$, consisting of the forbidden factors for $F$. 

Let
$W = \{{\bm f}\in B^\omega
: 01\bullet{\bm f}\in G\}.$ 
To prove Theorem~\ref{fife01} it is enough to prove that $W = F$.
Let $L\subseteq\Sigma^\omega$
 and let $x\in\Sigma^*$. We define the (left) quotient $x^{-1}L$ by
$x^{-1}L= \{{\bm y} \in  \Sigma^\omega 
: x{\bm y} \in  L\}$.
The next lemma establishes several
identities concerning quotients of the set $W$. They are proved using (\ref{2.3 01}) and Lemma~\ref{allouche}. The identities demonstrate that $W$ is precisely the set of infinite labeled paths through the automaton ${A_{01}}$ given in
Figure~\ref{automaton}. These are just the labelled paths omitting factors in $I$, so that  $W = F$. Thus, proving Lemma~\ref{identities} establishes Theorem~\ref{fife01}.

\begin{lemma}\label{identities} The following identities hold:
\begin{enumerate}\item[(a)]$W = \alpha^{-1}W$;
\item[(b)]$(\beta\alpha\alpha)^{-1}W=(\beta\beta)^{-1}W=\beta^{-1}W$
\item[(c)]$(\beta\alpha\beta)^{-1}W=\gamma^{-1}W=(\beta\gamma)^{-1}W=(\gamma\gamma)^{-1}W$
\item[(d)]$(\beta\alpha)^{-1}W=(\gamma\alpha)^{-1}W$
\item[(e)]$(\beta\alpha\gamma)^{-1}W=(\gamma\beta)^{-1}W=\emptyset$
\end{enumerate}
\end{lemma}
Each set of identities corresponds to the state of ${A_{01}}$ with the same label as the identities. The non-accepting sink (e) is not shown in the figure.

\begin{figure}\label{automaton}\caption{`Fife' automaton ${A_{01}}$ for $G_{01}$}
\centering \begin{tikzpicture}[shorten >=1pt,node distance=100pt,auto]
  \tikzstyle{every state}=[fill={rgb:black,1;white,10}]

  \node[state,initial]   (a)                      {$a$};
  \node[state] (b) [above right of=a]  {$b$};
  \node[state] (c) [below right of=a] {$c$};
  \node[state]           (d) [below right of=b]     {$d$};

  \path[<-]
  (a)   edge [loop right]    node {$\alpha$} (a)
  (b) edge  node {$\beta$} (a)
       edge  [loop above] node {$\beta$} (b)
       edge  [bend right] node {$\alpha$} (d)
  (c) edge   node {$\gamma$} (a)
       edge [bend right]  node {$\beta$} (d)
       edge [loop below]    node {$\gamma$} (c)
       edge   node {$\gamma$} (b)
  (d) edge [bend right]  node {$\alpha$} (b)
       edge  [bend right] node {$\alpha$} (c)
 ;
\end{tikzpicture}
\end{figure}
\begin{proof}
Let $01\bullet{\bm f} = {\bm x}$ be good.

\noindent (a) \[\alpha{\bm f}\in W\iff 01\bullet\alpha{\bm f}\in G\iff(01\bullet\alpha)\delta(01)^{-1}\delta({\bm x})\in G\iff\]
\[\stkout{0100(0100)}^{-1}\delta( {\bm x})\in G\iff {\bm x}\in G\iff01\bullet{\bm x}\in G\iff {\bm f}\in W\]
so that $\alpha^{-1}W=W$.

\noindent (b) \[\beta{\bm f}\in W\iff 01\bullet\beta{\bm f}\in G\iff(01\bullet\beta)\delta(01)^{-1}\delta({\bm x})\in G\iff\]
\[01~\stkout{0100(0100)}^{-1}\delta({\bm x})\in G\iff \delta(0{\bm x})\in G\iff 0{\bm x}\in G\]

\[\beta\beta{\bm f}\in W\iff 01\bullet\beta\beta{\bm f}\in G\iff(01\bullet\beta\beta)\delta^2(01)^{-1}\delta^2({\bm x})\in G\iff\]
\[010100~\stkout{01000101(01000101)}^{-1}\delta^2({\bm x})\in G\iff\delta(0\delta(0{\bm x}))\in G\iff\] \[
\delta(0\delta(0{\bm x}))\in G
\iff 0\delta(0{\bm x})\in G
\iff \] \[
10{\bm x}\in G\mbox{ or }0{\bm x}\in G_{001}\iff 0{\bm x}\in G\]
Here we use the fact that 01 is a prefix of ${\bm x}$, so that  $10{\bm x}\not\in G$ and $ 0{\bm x}\in 001(0,1)^\omega$.

\[\beta\alpha\alpha{\bm f}\in W\iff 01\bullet\beta\alpha\alpha{\bm f}\in G\iff(01\bullet\beta\alpha\alpha)\delta^3(01)^{-1}\delta^3({\bm x})\in G\iff\]
\[01~\stkout{0100010101000100(0100010101000100)}^{-1}\delta^3({\bm x})\in G\iff\delta(0\delta^2({\bm x}))\in G\] \[
\iff 10\delta^2({\bm x})\in G\mbox{ or }0\delta^2({\bm x})\in G_{001}\iff0\delta^2({\bm x})\in G
 \] \[
\iff 1\delta({\bm x})\in G\mbox{ or }
\delta({\bm x})\in G_{001}\iff 1\delta({\bm x})\iff 0{\bm x}\in G\]
Here again note that 01 is a prefix of $\delta({\bm x})$ and $\delta^2({\bm x})$.
We have shown that $(\beta\alpha\alpha)^{-1}W=(\beta\beta)^{-1}W=\beta^{-1}W$, as desired.

\noindent (c) \[\gamma{\bm f}\in W\iff 01\bullet\gamma{\bm f}\in G\iff(01\bullet\gamma)\delta(01)^{-1}\delta({\bm x})\in G\iff\]
\[0101~\stkout{0100(0100)}^{-1}\delta({\bm x})\in G\iff \delta(00{\bm x})\in G\iff 00{\bm x}\in G\]

\[\beta\alpha\beta{\bm f}\in W\iff 01\bullet\beta\alpha\beta{\bm f}\in G\iff(01\bullet\beta\alpha\beta)\delta^3(01)^{-1}\delta^3({\bm x})\in G\iff\]
\[0101000101~\stkout{0100010101000100(0100010101000100)}^{-1}\delta^3({\bm x})\in G\iff\] \[\delta(0\delta^2({0\bm x}))\in G\iff
\iff 0\delta^2(0{\bm x})\in G\iff
1\delta(0{\bm x})\in G\mbox{ or }\delta(0{\bm x})\in G_{001}
 \] \[\iff
1\delta(0{\bm x})\in G\iff \delta(00{\bm x})\in G\iff 00{\bm x}\in G\]

\[\beta\gamma{\bm f}\in W\iff 01\bullet\beta\gamma{\bm f}\in G\iff(01\bullet\beta\gamma)\delta^2(01)^{-1}\delta^2({\bm x})\in G\iff\]
\[0101000100~\stkout{01000101(01000101)}^{-1}\delta^2({\bm x})\in G\iff\delta(0\delta(00{\bm x}))\in G\iff\] \[
0\delta(00{\bm x})\in G
\iff 
100{\bm x}\in G\mbox{ or }00{\bm x}\in G_{001}\iff 100{\bm x}\in G\]

\[\gamma\gamma{\bm f}\in W\iff 01\bullet\gamma\gamma{\bm f}\in G\iff(01\bullet\gamma\gamma)\delta^2(01)^{-1}\delta^2({\bm x})\in G\]
\[\iff010101000100~\stkout{01000101(01000101)}^{-1}\delta^2({\bm x})\in G\]  \[ \iff\delta^2(100{\bm x})\in G
\iff 100{\bm x}\in G\]

Finally, suppose $00{\bm x}\in G$. Then 0001 is a prefix of  $00{\bm x},$ and by Theorem~\ref{main theorem} we can write $00{\bm x}=\delta(10{\bm y})$, some  $10{\bm y}\in G$. Then
by Lemma~\ref{slide}, $010{\bm y}\in G$. Thus $\delta(010{\bm y})=0100{\bm x}\in G$, and in particular, $100{\bm x}\in G$. We conclude that 
$100{\bm x}\in G\iff 00{\bm x}\in G$, which gives the desired result.

\noindent (d)
\[\beta\alpha{\bm f}\in W\iff 01\bullet\beta\alpha{\bm f}\in G\iff(01\bullet\beta\alpha)\delta^2(01)^{-1}\delta^2({\bm x})\in G\iff\]
\[01~\stkout{01000101(01000101)}^{-1}\delta^2({\bm x})\in G\iff\delta(0\delta({\bm x}))\in G\iff\] \[
0\delta({\bm x})\in G
\iff 
1{\bm x}\in G\mbox{ or }{\bm x}\in G_{001}\iff 1{\bm x}\in G\]

\[\gamma\alpha{\bm f}\in W\iff 01\bullet\gamma\alpha{\bm f}\in G\iff(01\bullet\gamma\alpha)\delta^2(01)^{-1}\delta^2({\bm x})\in G\iff\]
\[0101~\stkout{01000101(01000101)}^{-1}\delta^2({\bm x})\in G\iff\delta(00\delta({\bm x}))\in G\]  \[
\iff\delta^2(1{\bm x})\in G\iff 1{\bm x}\in G\]

Thus $(\beta\alpha)^{-1}W=(\gamma\alpha)^{-1}W$.

\noindent (e)

\[\beta\alpha\gamma{\bm f}\in W\iff 01\bullet\beta\alpha\gamma{\bm f}\in G\iff(01\bullet\beta\alpha\gamma)\delta^3(01)^{-1}\delta^3({\bm x})\in G\iff\]
\[ 010100010101000101~\stkout{0100010101000100(0100010101000100)}^{-1}\delta^3({\bm x})\in G\]\[\iff\delta(0\delta^2(00{\bm x}))\in G\iff 0\delta^2(00{\bm x})\in G
\iff
 \] \[ 1\delta(00{\bm x})\in G\mbox{ or }\delta(00{\bm x})\in G_{001}
\iff 1\delta(00{\bm x})\in G\iff
000{\bm x}\in G\]
But $000{\bm x}$ has prefix 0000. Thus $(\beta\alpha\gamma)^{-1}W=\emptyset$.

\[\gamma\beta{\bm f}\in W\iff 01\bullet\gamma\beta{\bm f}\in G\iff(01\bullet\gamma\beta)\delta^2(01)^{-1}\delta^2({\bm x})\in G\iff\]
\[01010100~\stkout{01000101(01000101)}^{-1}\delta^2({\bm x})\in G\iff\delta^2(10{\bm x})\in G\iff 10{\bm x}\in G\]

But $10{\bm x}$ has prefix 1001, so $(\gamma\beta)^{-1}W=\emptyset$.

\end{proof}

  We have
\[G=G_{01}\cup G_{001}\cup G_{0001}\cup G_1.\]
We can write Theorem~\ref{fife01} as \[G_{01}=01\bullet W.\] 

 By Lemma~\ref{slide}, ${\bm x}\in G_1$ if and only if 
$0{\bm x}\in G_{01}$, so that \[G_{1}=0^{-1}G_{01}.\] Also, by Theorem~\ref{main theorem} and Theorem~\ref{f good}, \[G_{0001}=\delta(G_{1})\] since $G_{10}=G_1$.
If ${\bm x}\in G_{001}$, we can write ${\bm x}=0\delta(0{\bm y})$ for some ${\bm y}$. However, by Lemma~\ref{allouche},
\begin{eqnarray*}
&&0\delta(0{\bm y})\in G\\
&\iff&10{\bm y}\in G\mbox{ or }0{\bm y}\in G_{001}\\
&\iff&010{\bm y}\in G\mbox{ or }0{\bm y}\in G_{001}\mbox{ by Lemma~\ref{slide}}\\
&\iff&010{\bm y}\in G_{01}\mbox{ or }0{\bm y}\in G_{001}\\
&\iff&0{\bm y}\in (01)^{-1}G_{01}\mbox{ or }0{\bm y}\in G_{001}\\
&\iff&0\delta(0{\bm y})\in0\delta((01)^{-1}G_{01}\cup\in G_{001}),\mbox{ so that}
\end{eqnarray*}
\[G_{001}=0\delta((01)^{-1}G_{01}\cup\in G_{001}).\]

In summary,
\begin{Theorem}\label{full fife}
\begin{eqnarray*}
G_{01}&=&01\bullet W\\
G_{001}&=&0\delta((01)^{-1}G_{01}\cup\in G_{001})\\
G_{0001}&=&\delta(G_{1})\\
G_{1}&=&0^{-1}G_{01}
\end{eqnarray*}
\end{Theorem}

Since $I$ is a regular language one could give an enumeration of finite prefixes of $G_{01}$ (and hence $G_{0001}$, $G_1$) following the approach of Kobayashi \cite{kobayashi}. This would give a lower bound on good words of length $n$.
Obtaining an enumeration of all prefixes of $G$ would involve dealing with the 
recursion in the equation for $G_{001}$, and would give a better lower bound on the number of length $n$ good words.

It would be nice to remove the recursion in the equation for $G_{001}$. With ${\cal I}$, ${\cal F}$, $\hat{\alpha}$, $\hat{\beta}$, $\hat{\gamma}$, $\hat{\bullet}$, etc.,  corresponding to $ I$, $ F$, $\alpha$, $\beta$, $\gamma$, $\bullet$, etc., but for overlaps, Rampersad's formulation of Fife's theorem has the following (unrecursive) form:

\begin{Theorem}
Let ${\bm x}\in\{0, 1\}^\omega$.
\begin{enumerate}
\item If ${\bm x}$ begins with 01, then ${\bm x}$ is overlap-free if and only if ${\bm x}= 01\hat{\bullet} {\bm f}$ for some
${\bm f}\in {\cal F}$.
\item If ${\bm x}$ begins with 001, then ${\bm x}$ is overlap-free if and only if ${\bm x}= 01\hat{\bullet} {\bm f}$ for some
$\beta{\bm f}\in {\cal F}$.
\end{enumerate}
\end{Theorem}

\begin{remark}In the second part of this theorem the condition can be rewritten as saying that ${\bm x}\in 001(0,1)^\omega$ is not overlap-free, if and only if $\hat{\beta}{\bm f}$ has a factor in ${\cal I}$, where ${\bm f}$ describes the canonical decomposition of 
${\bm x}$ . In particular, if ${\bm f}$ has a factor in ${\cal I}$, then ${\bm x}$
 is not overlap-free. No result analogous to this seems possible for good words; we shall see that there are words ${\bm x}\in G_{001}$ with description ${\bm x}=001\bullet{\bm f}$ such that ${\bm f}$ has a factor $\beta\alpha\beta\beta\in I$. On the other hand, $\alpha\alpha\alpha\gamma\not\in I$, but cannot be a factor of any ${\bm f}$ such that $001\bullet{\bm f}\in G.$
\end{remark}

Let ${\bm f} \in B^\omega$
such that $001 \bullet {\bm f} = {\bm x}$. Let $w\in B^k$. Then 

\begin{equation}\label{2.3 001}
001\bullet w{\bm f} = (001\bullet w)\delta^k(01)^{-1}\delta^k({\bm x})\end{equation}

Let $F_{001}$ be the set of infinite words walkable on the automaton ${A_{001}}$ of Figure~\ref{automaton2}. 

\begin{Theorem}\label{fife001}
Let ${\bm x}\in\{0, 1\}^\omega$.
If ${\bm x}$ begins with 001, then ${\bm x}$ is good if and only if ${\bm x}= 001\bullet f$ for some
${\bm f}\in F_{001}$.
\end{Theorem}
Let
$W_{001} = \{{\bm f}\in B^\omega
: 001\bullet{\bm f}\in G\}.$
We prove that $W_{001} = F_{001}$.

The  identities in the following lemma correspond to the states of ${A_{001}}$, except for state $\beta\gamma\alpha$, which is labelled for its shortest path from state a. The non-accepting sink (d) is omitted from the figure. Proving the identities thus proves 
Theorem~\ref{fife001}. One notes that $\beta\alpha\beta\beta\in I\cap F_{001}$. However, $\alpha\alpha\alpha\gamma$ is a prefix of words in $F$, but cannot be a factor of any ${\bm f}\in F_{001}$.

\begin{lemma}\label{identities} The following identities hold:
\begin{enumerate}\item[(a)]$(\beta\alpha)^{-1}W_{001}=W_{001} = \alpha^{-1}W_{001}$;
\item[(b)]$(\beta\gamma\alpha\alpha)^{-1}W_{001}=\beta\gamma\beta)^{-1}W_{001}=(\beta\beta)^{-1}W_{001}=\beta^{-1}W_{001}$
\item[(c)]$(\beta\gamma\alpha\beta)^{-1}W_{001}=\beta\gamma\gamma^{-1}W_{001}=(\beta\gamma)^{-1}W_{001}$
\item[(d)]$(\beta\gamma\alpha\gamma)^{-1}W_{001}=\gamma^{-1}=\emptyset$
\end{enumerate}
\end{lemma}
\begin{figure}\label{automaton2}\caption{`Fife' automaton ${A_{001}}$ for $G_{001}$}
\centering \begin{tikzpicture}[shorten >=1pt,node distance=100pt,auto]
  \tikzstyle{every state}=[fill={rgb:black,1;white,10}]

  \node[state,initial]   (a)                      {$a$};
  \node[state] (b) [right of=a]  {$b$};
  \node[state] (c) [above right of=b] {$c$};
  \node[state]           (d) [below right of=b]     {$\beta\gamma\alpha$};

  \path[->]
  (a)   edge [loop above]    node {$\alpha$} (a)
     edge     [bend right] node {$\beta$} (b)

  (b) edge  [loop below] node {$\beta$} (b)
edge  [bend right] node {$\alpha$} (a)       
edge  [bend right=5] node {$\gamma$} (c)
  (c) edge   node {$\alpha$} (d)
       edge [bend right=45]  node {$\beta$} (b)
       edge [loop above]    node {$\gamma$} (c)
  (d) edge [bend left=5]  node {$\alpha$} (b)
       edge  [bend right] node {$\beta$} (c)
 ;
\end{tikzpicture}
\end{figure}
\begin{proof}
Let $001\bullet{\bm f} = {\bm x}$ be good.

\noindent (a) \[\alpha{\bm f}\in W_1\iff 001\bullet\alpha{\bm f}\in G\iff(001\bullet\alpha)\delta(01)^{-1}\delta({\bm x})\in G\iff\]
\[0~\stkout{0100(0100)}^{-1}\delta( {\bm x})\in G\iff 1{\bm x}\in G\mbox{ or }{\bm x}\in G_{001}\iff{\bm x}\in G\iff {\bm f}\in W_{001}\]

\[\beta\alpha{\bm f}\in W_{001}\iff 001\bullet\beta\alpha{\bm f}\in G\iff(001\bullet\beta\alpha)\delta^2(01)^{-1}\delta^2({\bm x})\in G\iff\]
\[001~\stkout{01000101(01000101)}^{-1}\delta^2({\bm x})\in G\iff0\delta(0\delta({\bm x}))\in G\iff\] \[
10\delta({\bm x})\in G\mbox{ or }0\delta({\bm x})\in G_{001}
\iff 0\delta({\bm x})\in G
\iff 1{\bm x}\in G\mbox{ or }{\bm x}\in G_{001}\]\[\iff {\bm x}\in G\iff {\bm f}\in G\]

so that $(\beta\alpha)^{-1}W_{001}=\alpha^{-1}W_{001}=W_{001}$.

\noindent (b) \[\beta{\bm f}\in W_{001}\iff 001\bullet\beta{\bm f}\in G\iff(001\bullet\beta)\delta(01)^{-1}\delta({\bm x})\in G\iff\]
\[001~\stkout{0100(0100)}^{-1}\delta({\bm x})\in G\iff 0\delta(0{\bm x})\in G\iff 10{\bm x}\in G\mbox{ or }0{\bm x}\in G_{001}\iff 0{\bm x}\in G\iff 10{\bm x}\in G
\]

\[\beta\beta{\bm f}\in W_{001}\iff 001\bullet\beta\beta{\bm f}\in G\iff(01\bullet\beta\beta)\delta^2(01)^{-1}\delta^2({\bm x})\in G\iff\]
\[010100~\stkout{01000101(01000101)}^{-1}\delta^2({\bm x})\in G\iff\delta(0\delta(0{\bm x}))\in G\] \[
\iff 0\delta(0{\bm x})\in G
\iff 10{\bm x}\in G\mbox{ or }0{\bm x}\in G_{001}
\iff 
10{\bm x}\in G\]

\[\beta\gamma\alpha\alpha{\bm f}\in W_{001}\iff 001\bullet\beta\gamma\alpha\alpha{\bm f}\in G\iff(01\bullet\beta\gamma\alpha\alpha)\delta^4(01)^{-1}\delta^4({\bm x})\in G\iff\]
\[001 0100 0100 \delta^4({\bm x})\in G\]\[\iff0\delta(0\delta^2(1\delta({\bm x})))\in G \iff 10\delta^2(1\delta({\bm x}))\in G\mbox{ or }0\delta^2(1\delta({\bm x}))\in G_{001}
\iff\]\[ 0\delta^2(1\delta({\bm x}))\in G
\iff 1\delta(1\delta({\bm x}))\in G\mbox{ or }\delta(1\delta({\bm x})\in G_{001}\iff 1\delta(1\delta({\bm x}))\in G\]\[\iff 01\delta({\bm x})\in G\iff \delta(0{\bm x})\in G \iff 0{\bm x}\in G\iff 10{\bm x}\in G
\]

Thus $((\beta\gamma\alpha\alpha)^{-1}W_{001}=\beta\gamma\beta)^{-1}W_{001}=(\beta\beta)^{-1}W_{001}=(\beta\alpha)^{-1}W_{001}=\beta^{-1}W_{001}$.

\noindent (c) \[\beta\gamma{\bm f}\in W_{001}\iff 001\bullet\beta\gamma{\bm f}\in G\iff(001\bullet\beta\gamma)\delta^2(01)^{-1}\delta^2({\bm x})\in G\iff\]
\[001 01000100 ~\stkout{01000101(01000101)}^{-1}\delta^2({\bm x})\in G\iff 0\delta(0\delta(00{\bm x}))\in G\iff\]\[ 10\delta(00{\bm x})\in G\mbox{ or }0\delta(00{\bm x})\in G_{001}
\iff 0\delta(00{\bm x})\in G\]\[\iff 100{\bm x}\in G\mbox{ or }00{\bm x}\in G_{001}\iff 100{\bm x}\in G\]

\[\beta\gamma\gamma{\bm f}\in W_{001}\iff 001\bullet\beta\gamma\gamma{\bm f}\in G\iff(001\bullet\beta\gamma\gamma)\delta^3(01)^{-1}\delta^3({\bm x})\in G\iff\]
\[001 01000100 01000101 01000101~\stkout{0100010101000100(0100010101000100)}^{-1}\delta^3({\bm x})\in G\iff\] 
\[0\delta(0\delta^2(100{\bm x}))\in G
\iff 0\delta^2(100{\bm x})\in G\mbox{ or }\delta^2(100{\bm x})\in G_{001}\iff\]\[ 0\delta^2(100{\bm x})\in G
1\delta^(100{\bm x})\in G\mbox{ or }\delta^(100{\bm x})\in G_{001}
 \] \[\iff
1\delta^(100{\bm x})\in G\iff 0100{\bm x}\in G\iff 100{\bm x}\in G\]

\[\beta\gamma\alpha\beta{\bm f}\in W_{001}\iff 001\bullet\beta\gamma\alpha\beta{\bm f}\in G\iff(001\bullet\beta\gamma\alpha\beta)\delta^4(01)^{-1}\delta^4({\bm x})\in G\iff\]
\[001 01000100 0100010101000100 \delta^4({\bm x})\in G\iff 0\delta(0\delta^2(1\delta(0{\bm x})))\in G\]\[ 
\iff 10\delta^2(1\delta(0{\bm x}))\in G\mbox{ or }0\delta^2(1\delta(0{\bm x}))\in G_{001}
\iff 0\delta^2(1\delta(0{\bm x}))\in G\]\[
\iff 1\delta(1\delta(0{\bm x}))\in G\mbox{ or }\delta(1\delta(0{\bm x}))\in G_{001}\iff 1\delta(1\delta(0{\bm x}))\in G\]\[\iff 01\delta(0{\bm x})\in G\iff \delta(00{\bm x})\in G\iff 00{\bm x}\in G \iff 100{\bm x}\in G
\]

\noindent (d)

\[\gamma{\bm f}\in W_{001}\iff 001\bullet\gamma{\bm f}\in G\iff(001\bullet\gamma)\delta(01)^{-1}\delta({\bm x})\in G\iff00101~\stkout{0100(0100)}^{-1}\delta({\bm x})\in G\]\[\iff0\delta(00{\bm x}))\in G\] But $00{\bm x}$ has prefix 0000. Thus $\gamma^{-1}W_{001}=\emptyset$.

\[\beta\gamma\alpha\gamma{\bm f}\in W_{001}\iff 001\bullet\beta\gamma\alpha\gamma{\bm f}\in G\iff(001\bullet\beta\gamma\alpha\gamma)\delta^4(01)^{-1}\delta^4({\bm x})\in G\]
\[\iff 001 0100 0100 0100010101000100 0100010101000100 \delta^4({\bm x})\in G\iff\] \[0 \delta(0\delta(^2(1\delta(00{\bm x})))\in G
\]

We see that  $001\bullet\beta\gamma\alpha\gamma{\bm f} $contains the prefix 0000 of $00{\bm x}$, so that $(\beta\gamma\alpha\gamma)^{-1}W_{001}=\emptyset$.
\end{proof}

\section{Lexicographically extremal good words}
If $u$ is a word of positive length denote by $u^{-}$ the word obtained from $u$ by deleting its last letter. The lexicographic order on binary words is given recursively by
\[u< v\iff v\ne \epsilon\mbox{ and }((u=\epsilon)\mbox{ or }(u^-<v^-)\mbox{ or }((u=u^-0)\mbox{ and }(v=u^-1))).\]
Note that the morphism $\delta$ is order-reversing: Let $u$ and $v$ be non-empty binary words so that $u<v$. Write $u=u'0u^{\prime\prime}$, $v=u'1v^{\prime\prime}$ where $u'$ is the longest common prefix of $u$ and $v$. Then $\delta(u')01$ is a prefix of $\delta(u)$, while $\delta(u')00$ is a prefix of $\delta(v)$, so that $\delta(u)>\delta(v)$.

For each non-negative integer $n$, let $\ell_n$ (resp., $m_n$) be the lexicographically least (resp., greatest) word of length $n$ such that $\ell_n$ (resp., $m_n$) is the prefix of a one-sided infinite good word. 

\begin{lemma}\label{ell_n prefix} Let $n$ be a non-negative integer. Word $\ell_n$ is a prefix of $\ell_{n+1}$. Word $m_n$ is a prefix of $m_{n+1}$. 
\end{lemma}
\begin{proof} We prove the result for the $\ell_n$; the proof for the $m_n$ is similar. Let $\ell_n{\bf r}$ be a one-sided infinite good word. Let $p$ be the length $n+1$ prefix of
$\ell_n{\bf r}$, and let $q$ be the length $n$ prefix of
$\ell_{n+1}$. We need to show that $q=\ell_n$. Both $p$ and $q$ are prefixes of one-sided infinite good words. By definition we have $\ell_{n+1}\le p$ and $\ell_n\le q$. If $\ell_n<q$, then $p^-=\ell_n<q=\ell_{n+1}^-$, so that $p<\ell_{n+1}$. This is a contradiction. Therefore $\ell_n= q$, as desired.
\end{proof}

Let ${\bf {\bm\ell}}=lim_{n\rightarrow\infty}\ell_n$, ${\bf m}=lim_{n\rightarrow\infty}m_n$.

\begin{lemma} Word ${\bf {\bm\ell}}$ is the lexicographically least one-sided infinite good word. Word ${\bf m}$ is the lexicographically greatest one-sided infinite good word. 
\end{lemma}
\begin{proof} We show that ${\bf {\bm\ell}}$ is lexicographically least. The proof that ${\bf m}$ is lexicographically greatest is similar. Let ${\bf w}$ be a one-sided infinite good word. For each $n$ let $w_n$ be the length $n$ prefix of ${\bf w}$, so that ${\bf w}=lim_{n\rightarrow\infty}w_n$. 

If for some $n$ we have $w_n>\ell_n$, then ${\bf w}>{\bf {\bm\ell}}$. 

Otherwise $w_n\le \ell_n$ for all $n$. By the definition of the $\ell_n$ we have  $w_n\ge \ell_n$, so that $w_n= \ell_n$ for all $n$. Thus ${\bf w}=\lim_{n\rightarrow\infty}w_n=\lim_{n\rightarrow\infty}\ell_n={\bf {\bm\ell}}$. 

In all cases we find ${\bf w}\ge{\bf {\bm\ell}}$.
\end{proof}

\begin{lemma}\label{least} We have ${\bm\ell}=\delta({\bf m})$. 
\end{lemma}
\begin{proof} Word 0001 is the least good word of length 4. Also,  0001 is a factor of ${\bm d}$, so that there are one-sided infinite good words (which are suffixes of ${\bm d}$) with prefix 0001. Therefore, $\ell_4=0001$.
Since $0001$ is a prefix of ${\bm\ell}$, and ${\bm\ell}$ is good, we can write ${\bm\ell}=\delta({\bf m'})$ for some ${\bf m'}$ by Theorem~\ref{main theorem}. Since ${\bm\ell}$ is good, ${\bf m'}$ is good by Lemma~\ref{delta^-1 good}. It follows that ${\bf m'}\le {\bf m}$. However if ${\bf m'}< {\bf m}$ then $\delta({\bf m})<\delta({\bf m'})={\bm\ell}$ since $\delta$ is order-reversing. This is impossible, since ${\bm\ell}$ is least. Therefore ${\bf m'}= {\bf m}$, and ${\bm\ell}=\delta({\bf m})$.
\end{proof}

\begin{lemma}\label{greatest} We have $0{\bf m}=\delta({\bm\ell})$.
\end{lemma}
\begin{proof} Since 1 is the greatest good word of length 1, and ${\bm d}$ has suffixes which begin with 1, we have $m_1=1$. By Lemma~\ref{slide}, $0{\bf m}$ is good. Since 01 is a prefix of $0{\bf m}$, by Theorem~\ref{main theorem} we can write $0{\bf m}=\delta({\bm \ell'})$ for some ${\bm \ell'}$. Since $0{\bf m}$ is good, ${\bm \ell'}$ is good, so that ${\bm \ell'}\ge{\bm \ell}$. However, $\delta$ is order preserving, so that if 
${\bm \ell'}>{\bm \ell}$, then $\delta({\bm \ell})>\delta({\bm \ell'})={\bm m}$, contradicting the maximality of ${\bm m}$. Thus ${\bm \ell'}={\bm \ell}$ and $0{\bf m}=\delta({\bm \ell})$.
\end{proof}

\begin{Theorem}\label{lex least} Word ${\bm m}$ is the fixed point
\begin{equation}\label{m}{\bm m}=h_1({\bm m})\end{equation}
where $h_1=[1000,1010]$.
Word ${\bm \ell}$ is the fixed point 
\begin{equation}{\bm \ell}=h_2({\bm \ell}).\end{equation}
where $h_2=[0001,0101]$.
\end{Theorem}
\begin{proof}
From Lemma~\ref{least} and Lemma~\ref{greatest} we 
find ${\bm m}=0^{-1}\delta^2({\bm m})$, so that for each $n$ we have
$m_{4n-2}=0^{-1}\delta^2(m_n)$. 
However, $\delta^{2}(a)=0h_1(a)0^{-1}$ for $a\in\{0,1\}$, so that $\delta^2(w)=0h_1(w)0^{-1}$ for $w\in\{0,1\}^*$. 
Thus $m_{4n-2}=h_1(m_n)0^{-1}$.
Then \begin{eqnarray*}
{\bm m}&=&\lim_{n\rightarrow\infty} m_n\\
&=&\lim_{n\rightarrow\infty} m_{4n-2}\\
&=&\lim_{n\rightarrow\infty} h_1(m_n)0^{-1}\\
&=&\lim_{n\rightarrow\infty} h_1(m_n)\\
&=&h_1(\lim_{n\rightarrow\infty} m_n)\\
&=&h_1({\bm m})
\end{eqnarray*} 
The proof for ${\bm \ell}$ is similar, using the fact that $h_2=(01)^{-1}\delta^201$.
\end{proof}
\begin{corollary}Every finite factor of ${\bm \ell}$ or ${\bm m}$ is a factor of ${\bm d}$, and vice versa. However, words ${\bm \ell}$, ${\bm m}$, and ${\bm d}$ have no common suffix.
\end{corollary}
\begin{proof} Because $\delta^2$, $h_1$ and $h_2$ are conjugates of each other, their fixed points have the same factors. For example, let $u$ be a factor of ${\bm d}=\lim_{n\rightarrow\infty}\delta^{n}(0)$. Then for some $n$, $u$ is a factor of $\delta^{2n}(0)=0h_1^{2n}(0)0^{-1}$, which is a factor of $h_1^{2n}(1010)=
h_1^{2n+1}(1)$, which is a prefix of ${\bm m}=\lim_{n\rightarrow\infty}h_1^{n}(1)$. Thus $u$ is a factor of ${\bm m}$.

Now suppose ${\bm d}$ and ${\bm \ell}$ have a common suffix ${\bm s}$. (The proofs for the other pairs of fixed points are similar.) Write 
${\bm d}=d_i{\bm s}$ and ${\bm \ell}=\ell_j{\bm s}$ for some $i$ and $j$, where $d_n$ is the prefix of ${\bm d}$ of length $n$. Then 
\begin{eqnarray*}
d_i{\bm s}&=&{\bm d}\\
&=&\delta^2({\bm d})\\
&=&\delta^2(d_i)\delta^2({\bm s})
\end{eqnarray*} 
so that  ${\bm s}=s_{3i}\delta^2({\bm s})$ where is $s_n$ is the prefix of ${\bm s}$ of length $n$. Similarly,
\begin{eqnarray*}
\ell_j{\bm s}&=&{\bm \ell}\\
&=&h_2({\bm \ell})\\
&=&h_2(\ell_j)h_2({\bm s})
\end{eqnarray*} 
so that  ${\bm s}=s_{3j}h_2({\bm s}).$ Thus
\begin{eqnarray*}
s_{3j}h_2({\bm s})&=&{\bm s}\\
&=&s_{3i}\delta^2({\bm s})\\
&=&s_{3i}01h_2({\bm s})\\
&=&s_{3i+2}h_2({\bm s})\\
\end{eqnarray*}
It follows that $h_2({\bm s})$ has period $|3(i-j)+2|$. However, $h_2({\bm s})$ is good, and does not contain fourth powers. Therefore, $3(i-j)+2=0$. This is impossible, since $i$ and $j$ are integers.
\end{proof}
\section{Binary patterns in ${\bm d}$}

Every word encounters its factors as patterns. If a word encounters some binary pattern $p$, it necessarily encounters the complement $\overline{p}$. Shur \cite{shur96} has shown that, up to complementation, the only binary patterns encountered by ${\bf t}$ are its factors and 00100.  The situation with ${\bm d}$ is more complicated. 

\begin{lemma}\label{01 parse} Any factor $0u$ of ${\bm d}$ can be written as $h_1(p)$ for some word $p$ where $h_1=[0,01]$. Any factor $u0$ of ${\bm d}$ can be written as $h_2(p)$ for some word $p$ where $h_1=[0,10]$. 
Word ${\bm d}$ thus has an inverse image under each of $h_1$ and $h_2$.
\end{lemma}
\begin{proof} Every occurrence of 1 in ${\bm d}$ is preceded and followed by 0.
\end{proof}

\begin{corollary}\label{infinitely many} Let $0u$ be a factor of ${\bm d}$ such that $|0u|\ge 13$. Then $0u$ can be written as $h_1(p)$ where $p$ is neither a factor of ${\bm d}$ nor the complement of a factor of ${\bm d}$. Thus ${\bm d}$ encounters infinitely many patterns $p$ such that neither of $p$ and $\overline{p}$ is a factor of ${\bm d}$. 
\end{corollary}
\begin{proof}
The longest factor of ${\bm d}$ not containing $010001$ is 100010101000, which has length 12. Therefore, $0u=h_1(p)$ has the factor $010001=h_1(1001)$, and $p$ has the factor 1001. Neither of 1001 and 0110 is a factor of ${\bm d}$ so that neither $p$ nor $\overline{p}$ is a factor of ${\bf d}$.
\end{proof}

For any particular pattern $p$, the automatic proving system Walnut \cite{mousavi16} can in theory, given enough computing power and time, determine whether ${\bm d}$ encounters $p$. However, in the next theorem, we effectively characterize all binary patterns $p$ encountered by ${\bm d}$. The remainder of this section is devoted to its proof.

\begin{Theorem} Word $p$ is a binary pattern encountered by ${\bm d}$ if and only one of the following holds
\begin{enumerate}
\item One of $p$ and $\overline{p}$ is a factor of ${\bm d}$, $h_1({\bm d})$, or $h_2({\bm d})$.
\item One of $p$ and $\overline{p}$ is among
 
0010100,
01001001000,
00100100100,
001001001000,
00010010010,\\
000100100100,
0010001000100,
00100010001000, 
00010001000100,\\ and
000100010001000.
\end{enumerate}
The two possibilities are distinct.
\end{Theorem}

The following analogue of Lemma~\ref{uvu} will be useful for analyzing the patterns appearing in ${\bf w}$:
\begin{lemma}\label{uvu4} Let $u,v$ be binary words such that  $uvu$ is a factor of ${\bm d}$, $|u|\ge 3$, and either
\begin{itemize}
\item $|u|_{00}\ge 1$, or
\item $|u|_{10101}\ge 1$.
\end{itemize}
Then $|uv|\equiv 0$ (mod 4).
\end{lemma}
\begin{proof} Suppose that 00 is a factor of $u$, and $|u|\ge 3$. Then word $u$ contains one of the factors 000, 001, and 100. These can only arise in ${\bm d}$ as suffixes of some prefix of ${\bm d}$ of the form $\delta^2(p0)0$, $\delta^2(p0)01$, and $\delta^2(p0)$, respectively. The index of any occurrence of a factor 000, 001, or 100, and thus the index of any occurrence of $u$, in ${\bm d}$ is therefore fixed modulo 4, and the result follows.

Suppose $|u|_{10101}\ge 1$. The factor 10101 only occurs in ${\bm d}$ as a suffix of some prefix of ${\bm d}$ of the form $\delta^2(p10)(00)^{-1}$ and the index of any occurrence of $u$ in ${\bm d}$ is again fixed, modulo 4.
\end{proof}

\begin{Theorem}\label{patterns} Suppose that ${\bm d}$ encounters binary pattern $p$. Then one of $p$ and $\overline{p}$ either
\begin{enumerate}
\item  is a factor of ${\bm d}$, $h_1^{-1}({\bm d})$, or of $h_2^{-1}({\bm d})$, or
\item has the property that all its factors of the form $10^k1$ have the same length.
\end{enumerate}
\end{Theorem}
\begin{proof}
Without loss of generality, replacing $p$ by $\overline{p}$ if necessary, assume that 0 is the first letter of $p$. Assume that neither $p$ nor $\overline{p}$ is a factor of ${\bm d}$, $h_1({\bm d})$, or $h_2({\bm d})$. Since 000 is a factor of ${\bm d}$ and ${\bm d}$ doesn't encounter 0000,  we must have $|p|_1>0$.
Let $g(p)$ be a factor of $\delta^n(0)$ where $g=[X,Y]$ is a non-erasing morphism, and $n$ is as small as possible. By Lemma~\ref{even lengths}, $|X|$ and $|Y|$ are not both even. We make cases based on $|X|_1$, $|Y|_1$.

\begin{enumerate}
\item[{\bf Case 1}] We have $|X|_1$, $|Y|_1>0$:
If $XX$ (resp., $YY$, $XYX$, $YXY$) is a factor of $\delta^n(0)$, then by Lemma~\ref{uvu}, $|X|$ (resp., $|Y|$, $|XY|$, $|YXY|$) is even. 
It follows that not both $XX$ and $YY$ are factors of $g(p)$. Also, not both $XX$ and $XYX$ are factors of $g(p)$, or else $|X|$ and $|XY|$ are even, forcing $|Y|$ to be even. Similarly, not both $YY$ and $YXY$ are factors of $g(p)$.

If neither of $XX$ and $YY$ is a factor of $g(p)$, then $g(p)$ is an alternating string of $X$'s and $Y$'s, and thus a prefix of $XYXYXYX$. (Since ${\bm d}$ is good, the fourth power $(XY)^4$ is not a factor of ${\bm d}$.) But then $p$ is a prefix of 0101010, which is a factor of ${\bm d}$.  This is a contradiction.

Suppose then, that $XX$ is a factor of $g(p)$. Then neither of $YY$ and $XYX$ is a factor of $g(p)$. If $|p|_1\ge 2$, then $g(p)$ would have a prefix of the form $X^rYX^sY$, $r,s\ge 1$. This contains a factor $XYX$, which is impossible. It follows that $|p|_1\le 1$. Because ${\bm d}$ is good, $g(p)$ cannot have a factor $XXXX$, and we conclude that $p$ is a factor of 0001000. But 0001000 is seen to be a factor of ${\bm d}$. 

Finally, suppose that $YY$ is a factor of $g(p)$. Then $XX$ and $YXY$ are not factors of $g(p)$. Therefore, $p$ begins with 0, and has a factor 11, but not  00, 101, or 1111. (A factor 1111 in $p$ would give a fourth power in the good word ${\bm d}$.) It follows that $p$ is one of 011, 0110, 0111, and 01110. However $h_1(\overline{011})=0100$, $h_1(\overline{0110})=010001$, $\overline{0111}=1000$, and  $\overline{0111}=1000$ are all factors of ${\bm d}$.
\item [{\bf Case 2}]
We have $|X|_1=|Y|_1=0$. This forces $g(p)$ to be a factor of 000, so that $p$ is a binary word of length 3 or less. Each such word, or its complement, is a factor of ${\bm d}$. 
\item [{\bf Case 3 }] We have $|X|_1>0$ but $|Y|_1=0$. Write $Y=0^n$ where $1\le n\le 3$.
   \begin{enumerate}
   \item [{\bf Case 3(a)}] 
We have $|X|\ge 3$, and either $|X|_{00}\ge 1$ or
$|X|_{10101}\ge 1$.
If $p$ contains only a single $0$ (its first letter) then, since fourth powers do not appear in ${\bm d}$, $p$ must be a prefix of $0111$, which is the complement of a factor of ${\bm d}$. Otherwise, $p$ has a factor of the form $01^k0$ for some $k$, $0\le k\le 3$. 
Thus $XY^kX=X0^{nk}X$ is a factor of ${\bm d}$. By Lemma~\ref{uvu4} with $a=X$, $b=0^{nk}$, we have
$|X|\equiv -nk$ (mod 4). 

If $n=2$, then since $0^4$ is not a factor of ${\bm d}$, we have $0\le k\le 1$. Also, $|X|\equiv -nk$ (mod 4), so $|X|$ is even, giving $k\equiv -|X|/2$ (mod  2) and $k$ is determined; each pair of 0's in $p$ is separated by the same number of 1's. Then $\overline{p}$ has the property that all its factors of the form $10^k1$ have the same length.
If $n=1$ or $n= 3$, then since $0\le k\le 3$, $k$ is determined by the congruence $|X|\equiv -nk$ (mod 4), and again $\overline{p}$ has the property that all its factors of the form $10^k1$ have the same length.

   \item [{\bf Case 3(b)}] Either $|X|< 3$, or $|X|_{00}=|X|_{10101}=0$. If $|X|_{00}=0$, then $X$ doesn't contain a factor 00 or 11, and is therefore an alternating string of 0's and 1's. Thus the given conditions imply that $X=00$, or $X$ is a factor of 01010.

If $X=00$, then $g(p)$ must be a factor of 000, which is a factor of ${\bm d}$.

Suppose then that $X$ is a factor of 01010.
 If $p$ contained at most a single $0$ (its first letter), then $p$ would be a prefix of $0111$, the complement of a factor of ${\bm d}$. This is impossible. Therefore, assume that $p$ has a factor of the form $01^k0$. If $|X|\ge 3$, then $X$ is one of 010, 101, 0101, 1010, or 01010. One checks that for each of these possibilities for $X$ there is exactly one value of $k$ such that $X0^kX$ is a factor of ${\bm d}$:
\begin{itemize}
\item If $X$ is 010 or 01010, $k=1$.
\item If $X=101$, then a factor $XX$ would imply ${\bm d}$ has factor 11; a factor $X0X$ in ${\bm d}$ would extend on the right to a factor $X0X0=(10)^4$ in ${\bm d}$; a factor $X00X$ would imply ${\bm d}$ has factor 1001; thus $k=3$.
\item If $X$ is 0101 or 1010, then $XX$ is a fourth power, while $X0X$ has the factor 1001; $X000X$ contains 0000; thus $k=2$.
\end{itemize}
Thus if $|X|\ge 3$ then $\overline{p}$ has the property that all its factors of the form $10^k1$ have the same length.
 We are left with the cases where $|X|\le 2$ and $X\ne 00$, i.e., $X$ is among 1, 01, and 10.
\begin{itemize}
\item[$X=1$.]  Since ${\bm d}$ has a factor $XY^kX$, we cannot have $Y=00$ or else ${\bm d}$ contains 0000 or 1001. Thus $Y=0$. Then $g=[1,0]$, and $p$ is the complement of a factor of ${\bm d}$.  
\item[X=$01$.] Since $|X|$ is even, we must have $|Y|$ odd, so that $n=1$ or $n=3$. If $n=3$, then $YX=0^41$ and $YY=0^6$ are not factors of ${\bm d}$, forcing $|p|_1=1$, and $p$ is a factor of 00010000, which is a factor of ${\bm d}$. Suppose then that $n=1$. Now, however $g=[01,0]$, and $h_1(\overline{p})=g(p)$ is a factor of ${\bm d}$, so  that $\overline{p}$ is a factor of $h_1^{-1}({\bm d})$.
\item[X=$10$.] Again $|X|$ is even, so we find $n=3$ or $n=1$. In the first case, $XY$ contains 0000, which is impossible. Thus $p$ is a factor of $000$, a factor of ${\bm d}$. In the second case, $g=[10,0]$, and $\overline{p}$ is a factor of $h_2^{-1}({\bm d})$.
\end{itemize}
\end{enumerate}
\item [{\bf Case 4 }] We have $|X|_1=0$ but $|Y|_1>0$. Write $X=0^n$ where $1\le n\le 3$.
\begin{enumerate}
   \item [{\bf Case 4(a)}] 
We have $|Y|\ge 3$, and either $|Y|_{00}\ge 1$ or
$|Y|_{10101}\ge 1$.
If $p$ contains at most a single $1$ then, since fourth powers do not appear in ${\bm d}$, $p$ must be a factor of $0001000$, which is a factor of ${\bm d}$. Otherwise, $p$ has some factor of the form $10^k1$ for some $k$, $0\le k\le 3$. 
Thus $YX^kY=Y0^{nk}Y$ is a factor of ${\bm d}$. By Lemma~\ref{uvu4} with $a=Y$, $b=0^{nk}$, we have
$|Y|\equiv -nk$ (mod 4). 

If $n=2$, then since $0^4$ is not a factor of ${\bm d}$, we have $0\le k\le 1$. Also, $|Y|\equiv -nk$ so $|Y|$ is even, giving $k\equiv -|Y|/2$ (mod  2) and $k$ is determined. Therefore, $p$ has the property that all its factors of the form $10^k1$ have the same length.

If $n=1$ or $n= 3$, then since $0\le k\le 3$, $k$ is determined by the congruence $|Y|\equiv -nk$ (mod 4), and again $p$ has the property that all its factors of the form $10^k1$ have the same length.

\item [{\bf Case 4(b)}] Either $|Y|\le 3$, or $|Y|_{00}=|Y|_{10101}=0$. If $|Y|_{00}=0$, then $Y$ doesn't contain a factor 00 or 11, and is an alternating string of 0's and 1's. Thus the given conditions imply that $Y=00$, or $Y$ is a factor of 01010.

If $Y=00$, then $g(p)$ must be a factor of 000, which is a factor of ${\bm d}$.

Suppose then that $Y$ is a factor of 01010.
If $p$ contains at most a single $1$, then $p$ must be a factor of $0001000$, a factor of ${\bm d}$. Therefore, assume that $p$ has a factor of the form $10^k1$. If $|Y|>2$, an analysis as in Case 3(a),  shows that $p$ has the property that all its factors of the form $10^k1$ have the same length. The cases where $|Y|\le 2$ are impossible, with an analysis analogous to that of Case 3(a).
\end{enumerate}
\end{enumerate}
\end{proof}

\begin{lemma}\label{reversal} The set of factors of ${\bm d}$ is closed under reversal.
\end{lemma}
\begin{proof} To begin, we notice that if $u$ is a binary word, then
\[(\delta(u))^R=0^{-1}\delta(u^R)0.\]
Suppose $v$ is a factor of ${\bm d}$.  The set of length 2 factors of ${\bm d}$ is seen to be $\{00, 01, 10\}$, which is closed under reversal. Therefore, if $|v|\le 2$, then $v^R$ is a factor of ${\bm d}$. Suppose that $|v|>2$, and for every shorter factor $u$ of ${\bm d}$, the reversal $u^R$ is a factor of ${\bm d}$. However, $|v|>2$ implies that $v$ is a factor of $\delta(u)$, some factor $u$ of ${\bm d}$ which is shorter than $v$. Thus $u$ and $u^R$ are factors of ${\bm d}$. Then $(\delta(u))^R=0^{-1}\delta(u^R)0$ is also a factor of ${\bm d}$, so that $v^R$ is also. The result follows by induction.
\end{proof}
\begin{corollary}\label{reverse} If ${\bm d}$ encounters pattern $p$, then ${\bm d}$ encounters pattern $p^R$ also.
\end{corollary}

\begin{Theorem} The following are equivalent:
\begin{enumerate}
\item Word $p$ is a binary pattern encountered by ${\bm d}$, but neither $p$ nor $\overline{p}$ is a factor of ${\bm d}$, $h_1({\bm d})$, or $h_2({\bm d})$.
\item One of $p$ and $\overline{p}$ is among
 
0010100,
01001001000,
00100100100,
001001001000,
00010010010,\\
000100100100,
0010001000100,
00100010001000, 
00010001000100,\\ and
000100010001000.
\end{enumerate} 
\end{Theorem}
\begin{proof} $(1\Rightarrow 2)$: Suppose that $p$ is a binary pattern encountered by ${\bm d}$, but neither $p$ nor $\overline{p}$ is a factor of ${\bm d}$, $h_1({\bm d})$, or $h_2({\bm d})$. By Theorem~\ref{patterns}, replacing $p$ by $\overline{p}$ if necessary, there is a number $k$ such that there are exactly $k$ 1's between subsequent 0's in $p$. Thus $p=0^r(10^k)^s10^t$ for some $k,r,s,t$. Since ${\bm d}$ does not contain fourth powers, $k,r,s,t\le 3$. For the same reason, if $s=3$, then $r+s<k$. If $s=0$, then $p$ is a factor of 0001000, which is a factor of ${\bm d}$. We therefore have $s\ge 1$. We make cases based on the value of $k$:

\subsubsection*{$k=0$} In this case to avoid fourth powers we must also have $s\le 2$, and $p$ is a factor of 00011000 or 000111000. We note that  $h_1(\overline{00011000})=01010100010101$ 
is a factor of ${\bm d}$. Thus $p$ cannot be a factor of 00011000, which would make $\overline{p}$ a factor of $h_1^{-1}({\bm d})$. We conclude that $p$ is a factor of 000111000.
Note that, $h_2(00011100)=h_1(00111000)=00010101000$ is a factor of ${\bm d}$, so that $p$ cannot be a proper factor of 000111000. It remains that we must have $p=000111000$. To obtain a contradiction, we show that ${\bm d}$ does not encounter $p$. 

Suppose $g(p)$ is a factor of ${\bm d}$, $g=[X,Y]$ a non-erasing morphism. By Lemma~\ref{even lengths} assume that one of $|X|$ and $|Y|$ is odd. Both $XX$ and $YY$ are factors of $g(p)=XXXYYYXXX$, so by Lemma~\ref{uvu} we must have $|X|_1=0$ or $|Y|_1=0$. 

Suppose $|X|_1=0$. Thus $XXX$ ends in 000. Since $XXXY$ is a factor of $g(p)$, which is a factor of ${\bm d}$, the first letter of $Y$ is 1. Since $XXX$ begins with 000 and $YXXX$ is a factor of $g(p)$, the last letter of $Y$ is 1. But then $YY$ contains the factor 11, so that 11 is a factor of the good word ${\bm d}$, which is impossible.

Suppose $|Y|_1=0$. Switching $X$ and $Y$ in the argument of the previous paragraph shows that this is impossible also.

\subsubsection*{$k=1$} We make subcases based on whether $s$ is 1, 2, or 3:

\subsubsection*{~~$s=1$} If $r+t\ge 5$, then $p$ contains 00010100 or its reverse. However
${\bm d}$ is good, and therefore doesn't encounter 00010100. By Corollary~\ref{reverse}, ${\bm d}$ doesn't encounter the reverse of ${\bm d}$. We conclude that $r+t\le 4$, so that $p$ is a factor of 0101000, 0010100, or 0001010. However, 0101000 and 0001010 are factors of ${\bm d}$, so $p$ cannot be a factor of one of these. Therefore $p$ must be a factor of 0010100. Both 010100 and 010100 are factors of ${\bm d}$, so $p$ cannot be a proper factor of 0010100, forcing $p=0010100$.

\subsubsection*{~~$s=2$} In this case $p$ is a factor of 00010101000, which is a factor of ${\bm d}$. 

\subsubsection*{~~$s=3$} If $r=t=0$, then $p=1010101$ is the complement of a factor of ${\bm d}$. However, if $r>0$ or $t>0$, then $p$ contains one of the fourth powers 01010101 and 10101010, so that ${\bm d}$ doesn't encounter $p$, which is a contradiction.

\subsubsection*{$k=2$} We make subcases based on whether $s$ is 1, 2, or 3:

\subsubsection*{~~$s=1$} If $r+t\le 5$ then $p$ is factor of $000100100$ or of $001001000$. However $h_2(000100100)=00010001000=h_1(001001000)$ is a factor of ${\bm d}$, so $p$ is a factor of $h_1^{-1}({\bm d})$ or of $h_2^{-1}({\bm d})$, a contradiction. If $r=t=3$, then $p=0001001000=h_2(00010100)$, and ${\bm d}$ encounters 00010100. Since ${\bm d}$ is good, this is impossible.

\subsubsection*{~~$s=2$} If $r+t\le 3$, then $p$ is a factor of 1001001000, 0100100100, 0010010010, or 0001001001; however each of $h_1(1001001000)=0100010001000$, $h_2(0100100100)=0100010001000$, $h_1(0010010010)=0001000100010$, and $h_2(0001001001)=0001000100010$ is a factor of ${\bm d}$, so that $p$ is a factor of $h_1^{-1}({\bm d})$ or of $h_2^{-1}({\bm d})$. 

If $r=t=3$, then $p=0001001001000$. To get a contradiction, we show that ${\bm d}$ doesn't encounter $p$. Otherwise, suppose without loss of generality that ${\bm d}$ has a factor $g(p)$ where $g=[X,Y]$ is a non-erasing morphism and  $|X|$, $|Y|$ are not both odd. Word ${\bm d}$ has the factor $XXXYXXYXXYXXX$. We cannot have $|X|_1>0$, or since $XX$ and $XYX$ are factors of ${\bm d}$,  Lemma~\ref{uvu} would force both $|X|$ and $|Y|$ to be even. Write $X=0^n$. Since $XX$ is a factor of ${\bm d}$, we must have $n=1$ and $X=0$. Since $XXXY$ and $YXXX$ are factors of ${\bm d}$,  the first and last letters of $Y$ are 1. Then $YXXY$ contains the factor 1001, which cannt be a factor of ${\bm d}$. This is a contradiction.

We conclude that $4\le r+t\le 5$, so that $p$ is one of 01001001001000, 00100100100100, 00010010010010, 001001001001000, and 000100100100100.

\subsubsection*{~~$s=3$} If $r+t\le 1$, then $p$ is a factor of 01001001001 or 10010010010. However $h_2(01001001001)=010001000100010$ and $h_1(010001000100010)$ are factors of ${\bm d}$, so $p$ is a factor of $h_1^{-1}({\bm d})$ or of $h_2^{-1}({\bm d})$, which is a contradiction. On the other hand, if $r+t\ge 2$, then one of the fourth powers $(001)^4)$, $(010)^4$ and $(100)^4$ is a factor of $p$, and ${\bm d}$ contains a fourth power, which is impossible.

\subsubsection*{$k=3$} We make subcases based on whether $s$ is 1, 2, or 3:

\subsubsection*{~~$s=1$} In this case $p$ is a factor of 00010001000, which is a factor of ${\bm d}$. 

\subsubsection*{~~$s=2$} If $r\le 1$ then $p$ is a factor of 0100010001000 which is a factor of ${\bm d}$. This is impossible. Similarly $t\le 1$ is impossible. Thus $p$ must be one of 0010001000100, 00010001000100, 00100010001000, and 000100010001000

\subsubsection*{~~$s=3$} If $r,t\le 1$ then $p$ is a factor of 010001000100010, which is a factor of ${\bm d}$, a contradiction. Suppose $r\ge 2$ or $t\ge 2$. Then $p$ has $q=001000100010001$ or its reverse as a factor. By Corollary~\ref{reverse}, it suffices to show that ${\bm d}$ does not encounter $q$.

To get a contradiction, suppose that $g(q)$ is a factor of ${\bm d}$ for some non-erasing morphism $g=[X,Y]$ where one of $|X|$ and $Y$ is odd. Thus word ${\bm d}$ has the factor  $XXYXXXYXXXYXXXY$. Since $XX$ and $XYX$ are factors of ${\bm d}$, we must have $X=0^n$ for some $n$; otherwise, Lemma~\ref{uvu} implies that $|X|$ and $|Y|$ are even, a contradiction. Since $XX$ is a factor of ${\bm d}$, $n=1$ and $X=0$. Since $YXXX=Y000$ is a factor of ${\bm d}$, the last letter of $Y$ is 1. Therefore, $Y$ is always followed by 0 in ${\bm d}$, so $XXYXXXYXXXYXXXY0$ is a factor of ${\bm d}$. However $XXYXXXYXXXYXXXY0=(00Y0)^4$, and ${\bm d}$ contains a fourth power. This is impossible.

\noindent $(1\Leftarrow 2)$: We show that ${\bm d}$ encounters each of the listed patterns, but none of these patterns or their complements is a factor of ${\bm d}$, $h_1({\bm d})$, or $h_1({\bm d})$.

First we show that ${\bm d}$ encounters the listed patterns. Note that each of these patterns is   a factor of 
0010100, 
001001001000, 000100100100, or
000100010001000. Since 001001001000 is the reverse of 000100100100, by Lemma~\ref{reverse} it suffices to show that ${\bm d}$ encounters 
0010100, 
001001001000, or
000100010001000.
\begin{itemize}
\item $p=0010100$: Word ${\bm d}$ has the factor \[00010001000=g(0010100)\] where $g=[0,010]$.

\item $p=001001001000$: Word ${\bm d}$ has the factor \[0001010100010101000010101000=g(001001001000)\] where $g=[0,010101]$.

\item $p=000100010001000$: Word ${\bm d}$ has the factor \[0001010100010101000010101000=g(000100010001000)\] where $g=[0,10101]$.
\end{itemize}

Next we show that 
none of these patterns or their complements is a factor of ${\bm d}$, $h_1^{-1}({\bm d})$, or $h_2^{-1}({\bm d})$.
Notice that $h_1(101)=01001$ and $h_2(101)=10010$, and 1001 is not a factor of ${\bm d}$, since ${\bm d}$ is good. The complements of all the listed patterns contain 101 as a factor, so none of the complements is a factor of $h_1^{-1}({\bm d})$, or $h_2^{-1}({\bm d})$. Further, all the complements of the listed factors contain 11 as a factor, so none of the complements is a factor of ${\bm d}.$
It therefore suffices to show that none of the listed patterns is a factor of ${\bm d}$, $h_1^{-1}({\bm d})$, or $h_2^{-1}({\bm d})$.

Each of the listed patterns has one of 
0010100, 
01001001000,
00100100100,
00010010010,
 and
000100010001000 as a factor. 
It therefore suffices to show that none of these five patterns is a factor of ${\bm d}$, $h_1^{-1}({\bm d})$, or $h_2^{-1}({\bm d})$. 

\subsubsection*{Case 1: 0010100.} If this was a factor of ${\bm d}$, the 1's would force it to appear in the context $00010100$, which is impossible since ${\bm d}$ is good. Since 101 is a factor of 0010100, it is not a factor of $h_1^{-1}({\bm d})$, or $h_2^{-1}({\bm d})$ either.

\subsubsection*{Case 2: 01001001000.} This has 1001 as a factor, and therefore is not a factor of ${\bm d}$.

Word $00100010001000=h_1(01001001000)$ cannot be a factor of ${\bm d}$; if it was followed by 0 then 0000 would be a factor of ${\bm d}$; if it was followed by 10 then $(0100)^4$ is a factor of ${\bm d}$, which is impossible.

Word $h_2(01001001000)$ ends in the fourth power 0000 and cannot be a factor of ${\bm d}$

\subsubsection*{Case 3: 00100100100.} This has 1001 as a factor, and therefore is not a factor of ${\bm d}$.

Word $h_1(00100100100)=00010001000100$ cannot be a factor of ${\bm d}$; if it was followed by 1 then 1001 would be a factor of ${\bm d}$; if it was followed by 00 then $0000$ is a factor of ${\bm d}$; if it was followed by 01 then $(0001)^4$ is a factor of ${\bm d}$, which is impossible.

Word $00100010001000=h_2(00100100100)$ cannot be a factor of ${\bm d}$ as argued in Case 2.

\subsubsection*{Case 4: 00010010010.} This has 1001 as a factor, and therefore is not a factor of ${\bm d}$.

Word $h_1(00010010010)$ begins with 0000 and cannot be a factor of ${\bm d}$.

Word $00010001000100=h_2(00010010010)$ cannot be a factor of ${\bm d}$ as argued in Case 3.

\subsubsection*{Case 5: 000100010001000.} This cannot be a factor of ${\bm d}$; if followed by 0 it gives a factor 0000; if followed by 1 it gives a factor $(0001)^4$, which is impossible.

Word $h_1(000100010001000)$ begins with 0000 and cannot be a factor of ${\bm d}$.

Word $h_2(000100010001000)$ ends with 0000 and cannot be a factor of ${\bm d}$.
\end{proof}

\end{document}